\documentclass[a4paper,11pt]{article}

\usepackage{graphicx}
\usepackage{amssymb,amsmath,amsthm,mathrsfs,xfrac}
\usepackage{enumerate}
\usepackage{a4wide}
\usepackage{hyperref}
\usepackage{array}

\numberwithin{equation}{section}
\allowdisplaybreaks[4]

\newtheorem{proposition}{Proposition}[section]
\newtheorem{theorem}[proposition]{Theorem}
\newtheorem{lemma}[proposition]{Lemma}
\newtheorem{definition}[proposition]{Definition}

\newtheorem{remark}[proposition]{Remark}

\renewenvironment{proof}{\smallskip\noindent\emph{\textbf{Proof.}}%
  \hspace{1pt}}{\hspace{-5pt}{\nobreak\quad\nobreak\hfill\nobreak%
    $\square$\vspace{2pt}\par}\smallskip\goodbreak}

\newenvironment{proofof}[1]{\smallskip\noindent{\textbf{Proof~of~#1.}}%
  \hspace{1pt}}{\hspace{-5pt}{\nobreak\quad\nobreak\hfill\nobreak%
    $\square$\vspace{2pt}\par}\smallskip\goodbreak}

\newcommand{\pint}[1]{\mathaccent23{#1}}

\newcommand{\C}[1]{\mathbf{C}^{#1}}
\newcommand{\Cc}[1]{\mathbf{C}_c^{#1}}

\newcommand{\BV}{\mathbf{BV}}

\renewcommand{\L}[1]{{\mathbf{L}^#1}}
\newcommand{\Lloc}[1]{{\mathbf{L}_{\mathbf{loc}}^{#1}}}

\newcommand{\modulo}[1]{{\left|#1\right|}}
\newcommand{\norma}[1]{{\left\|#1\right\|}}
\newcommand{\caratt}[1]{{\chi_{\strut#1}}}
\newcommand{\reali}{{\mathbb{R}}}
\newcommand{\naturali}{{\mathbb{N}}}

\renewcommand{\epsilon}{\varepsilon}
\renewcommand{\phi}{\varphi}
\renewcommand{\theta}{\vartheta}

\newcommand{\sgn}{\mathop{\rm sgn}}
\newcommand{\sgnp}{\operatorname{sgn}^{+}}
\newcommand{\sgnm}{\operatorname{sgn}^{-}}

\renewcommand{\d}[1]{\mathinner{\mathrm{d}{#1}}}
\renewcommand{\div}{\mathinner{\mathop{{\rm div}}}}

\DeclareMathOperator*{\esssup}{ess\,sup}
\DeclareMathOperator*{\esslim}{ess\,lim}

\DeclareMathOperator*{\tr}{tr}

\makeatletter
\let\@fnsymbol\@arabic
\makeatother

\title{Definitions of solutions to the IBVP \\for multiD scalar balance
  laws}

\author{Elena Rossi\footnote{INdAM Unit, University of Brescia, Italy
    \texttt{elena.rossi@unibs.it}}}
\date {}

\begin{document}

\maketitle

\begin{abstract}
  \noindent
  We consider four definitions of solution to the initial--boundary
  value problem for a scalar balance laws in several space
  dimensions. These definitions are generalised to the same most
  general framework and then compared. The first aim of this paper is
  to detail differences and analogies among them. We focus then on the
  ways the boundary conditions are fulfilled according to each
  definition, providing also connections among these various
  modes. The main result is the proof of the equivalence among the
  presented definitions of solution.

  \medskip

  \noindent\textit{2010~Mathematics Subject Classification: 35L65,
    35L04, 35L60}

  \medskip

\end{abstract}

\section{Introduction}
\label{sec:Intro}

This paper is concerned with the relations among different definitions
of solution to the Initial--Boundary Value Problem (IBVP) for a
general scalar balance law in several space dimensions:
\begin{equation}
  \label{eq:1}
  \left\{
    \begin{array}{l@{\qquad}r@{\,}c@{\,}l}
      \partial_t u (t, x)
      +
      \nabla \cdot f\left(t,x,u (t,x)\right)
      =
      F \left(t,x,u (t,x)\right)
      & (t,x)
      & \in
      & \reali_+ \times \Omega
      \\
      u (0, x) = u_o (x)
      & x
      & \in
      & \Omega
      \\
      u (t,\xi) = u_b (t,\xi)
      & (t, \xi)
      & \in
      & \reali_+ \times \partial\Omega.
    \end{array}
  \right.
\end{equation}
Above and hereinafter, $\Omega$ is an open bounded subset of
$\reali^N$, with smooth boundary $\partial\Omega$, and
$\reali_+ = \left[0, +\infty\right[$.  The way the boundary condition
is satisfied is going to be precised further on and constitutes a key
issue addressed in this paper.

\smallskip

The pioneering work by Bardos, le Roux and
N\'ed\'elec~\cite{BardosLerouxNedelec} introduces a definition of
solution to~\eqref{eq:1} following the spirit of the one given by
Kru\v zkov in~\cite{Kruzkov} in the case without boundary. The idea of
the authors is to include in a unique integral inequality both Kru\v
zkov definition and the boundary condition. However, the BLN--definition
considers only functions admitting a trace at the boundary, for
instance $\BV$ functions. In~\cite{BardosLerouxNedelec}, the authors
explain the way the boundary condition has to be understood and
introduce a key inequality on the boundary, which we call the BLN condition,
relating the boundary datum to the trace of the solution,
see~\eqref{eq:blnbc}.

It is also possible to consider a definition of solution
to~\eqref{eq:1} analogous to the BLN--one, though involving classical
(regular) entropy--entropy flux pairs,
see~\cite[Definition~2.5]{bordo} and Definition~\ref{def:esol}. In
this way there is a sort of symmetry between, on one side, the
BLN--definition and this one for IBVPs as~\eqref{eq:1} and, on the
other side, Kru\v zkov definition and the definition of weak entropy
solution for initial value problems on all $\reali^N$.

On all $\reali^N$, solutions to the initial value problem for a
general scalar balance law are usually found in $\L\infty$, therefore
a question naturally arises: is it possible to find a concept of
solution to the IBVP~\eqref{eq:1} in this function space?  The key
difficulty is that, in general, a function in $\L\infty$ does not
necessarily admit a trace at the boundary. A first proposal to
overcome this issue is given by Otto in his PhD thesis (a summary is
presented in~\cite{Otto}, while more details and proofs can be found
in~\cite[Chapter~2]{MalekEtAlBook}). In the case of autonomous scalar
conservation laws, Otto replaces the Kru\v zkov entropy--entropy flux
pairs as exploited in~\cite{BardosLerouxNedelec} with the so-called
boundary entropy--entropy flux pairs and bases on them his definition
of solution. In this paper we provide a generalisation of this
definition to deal with non autonomous fluxes and arbitrary source
terms, see Definition~\ref{def:resol}.

Still looking for solutions to~\eqref{eq:1} in $\L\infty$, in the case
of scalar conservation laws with divergence free flux, Vovelle
introduces in~\cite{Vovelle} a definition of solution using the
so-called Kru\v zkov semi-entropy--entropy flux pairs. This definition
is then extended by Martin in~\cite{Martin} to deal with general
scalar balance laws. The resulting MV--definition, exploiting Kru\v
zkov semi entropies, resembles the BLN--definition, although not
requiring the existence of the trace of the solution at the boundary.

To sum up, there are mainly four definitions of solution
to~\eqref{eq:1}, which can be classified as follows: involving the
trace at the boundary (BLN--solutions and classical
entropy--solutions, i.e.~Definition~\ref{def:esol}) or not (Otto-type
solutions, i.e.~Definition~\ref{def:resol}, and MV--solutions);
dealing with regular entropies (entropy--solutions and Otto--type
ones) or with Lipschitz ones (BLN and MV--solutions). Definition~\ref{def:resol} and MV--definition share the interesting
feature of being stable under $\L1$-convergence, see
Remark~\ref{rem:stableL1} and also~\cite[Chapter~2,
Remark~7.33]{MalekEtAlBook}.

\smallskip

In this paper, we prove first the equivalence of
Definition~\ref{def:resol} and MV--definition, and then focus on the
way the boundary conditions are fulfilled according to those
definitions. We proceed similarly for the definitions of solution
requiring the existence of the trace, that is entropy--solutions and
BLN--solutions.

The main result of this paper is the proof of the equivalence among
all the presented definitions of solution to~\eqref{eq:1}. Of course,
this can be done only when the existence of the trace of the solution
at the boundary is assumed. Further information on the existence of
the trace at the boundary can be found
in~\cite{Panov2005,Panov2007,Vasseur}, see Section~\ref{sec:solTR} for
more details.
As an intermediate step, we also
prove the equivalence among the way the boundary conditions are
understood according to the various definitions.

\medskip

The paper is organised as follows. Section~\ref{sec:not} collects the
notation used throughout the paper. Sections~\ref{sec:def}
and~\ref{sec:bc} are devoted to the Otto-type definition of solution
and Martin--Vovelle--one: the first section contains the definitions
themselves and the theorem stating the equivalence between them, while
results on the way the boundary conditions are fulfilled constitute
the latter one. Section~\ref{sec:solTR} deals with the definitions of
solution with traces and provides also the main equivalence result.
In Section~\ref{sec:strongsol} we give the definition of strong
solution to~\eqref{eq:1} and a related results. Section~\ref{sec:1d}
provides further details on the one dimensional case, while
Section~\ref{sec:existence} summarises the existence results that can
be found in the literature. We collect the detailed proofs of our
results in Section~\ref{sec:TD}.




\section{Notation}
\label{sec:not}

The space dimension $N$, with $N \geq 1$, is fixed throughout. We set
$\reali_+=[0,+\infty[$. We denote by $\nu (\xi)$ the exterior normal
to $\xi \in \partial\Omega$. For $w,\, k \in \reali$ set
\begin{equation}
  \label{eq:I}
  \mathcal{I}[w,k] = \left\{
    z \in \reali \colon (w-z) (z-k) \geq 0
  \right\} =
  \left\{
    \theta \, w + (1 - \theta) k \colon \theta \in [0,1]
  \right\}.
\end{equation}
In other words, $\mathcal{I}[w,k]$ denotes the closed interval with
end points $w$ and $k$.

For the divergence of a vector field, possibly composed with another
function, we use the following notation:
\begin{displaymath}
  \nabla \cdot f\left(t,x,u (t,x)\right)
  = \div f\left(t,x,u (t,x)\right) +
  \partial_u f\left(t,x,u (t,x)\right) \cdot \nabla u (t,x).
\end{displaymath}
We use below the following standard assumptions:
\begin{description}
\item[(IC)] $u_o \in \L\infty (\Omega;\reali)$;
\item[(BC)] $u_b \in \L\infty ([0,T] \times \partial \Omega;\reali)$;
\item[(f)]
  $f \in \C2 ([0,T] \times \overline \Omega \times \reali; \reali^N)$,
  $\partial_u f \in \Lloc\infty ([0,T] \times \Omega \times \reali;
  \reali^N)$
  and
  $\nabla \cdot \partial_u f \in \L\infty ([0,T] \times \Omega \times
  \reali; \reali^{N\times N})$;
\item[(F)]
  $F \in \C2 ([0,T] \times \overline \Omega \times \reali; \reali)$
  and
  $\partial_u F \in \L\infty ([0,T] \times \Omega \times \reali;
  \reali)$.
\end{description}
Following~\cite{Martin, Vovelle}, we set
\begin{equation}
  \label{eq:sgn}
  \sgnp (s)
  =
  \begin{cases}
    1 & \mbox{if } s>0,
    \\
    0 & \mbox{if } s \leq 0,
  \end{cases}
  \qquad
  \sgnm (s)
  =
  \begin{cases}
    0 & \mbox{if } s \geq 0,
    \\
    -1 & \mbox{if } s < 0,
  \end{cases}
  \qquad
  \begin{array}{r@{\;}c@{\;}l}
    s^+
    & =
    & \max\{s, 0\} ,
    \\
    s^-
    & =
    & \max\{-s, 0\} .
  \end{array}
\end{equation}
We often use below the equalities $(-s)^- = s^+$ and $(-s)^+ = s^-$.
Introduce moreover the following notation: if
$g: \reali^2 \to \reali$, for all $z, \, w \in \reali$, set
\begin{align*}
  \partial_1 g (z,w) = \
  & \lim_{h \to 0}\frac{g (z+h,w)-g (z,w)}{h},
  &
    \partial_2 g (z,w) = \
  & \lim_{h \to 0}\frac{g (z,w+h)-g (z,w)}{h},
\end{align*}
and similarly for functions of more arguments.

\section{$\L1$-Stable Definitions}
\label{sec:def}

Before introducing the first definition of solution to~\eqref{eq:1},
we need to recall the notion of (classical) \emph{entropy--entropy
  flux pair}, see~\cite[Chapter~2, Definition~3.22]{MalekEtAlBook}.
\begin{definition}
  \label{def:eef}
  The pair
  $(\eta,q) \in \C2 (\reali;\reali) \times \C2 ([0,T]\times \overline
  \Omega \times \reali; \reali^N)$
  is called an \emph{entropy--entropy flux pair} with respect to $f$
  if
  \begin{enumerate}[i)]
  \item $\eta$ is convex, i.e.~$\eta'' (z)\geq 0$ for all
    $z\in \reali$;
  \item for all $t\in [0,T]$, for all $ x \in \Omega$, for all
    $x\in\reali$,
    $\partial_3 q (t,x,z) = \eta' (z) \, \partial_3 f (t,x,z) $.
  \end{enumerate}
\end{definition}
The notion of \emph{boundary entropy--entropy flux pair} is first
introduced by Otto in~\cite{Otto}, see also~\cite{MalekEtAlBook}, for
autonomous scalar conservation laws on bounded domains, and then
extended to a more general case in~\cite{Martin, Vovelle}. We recall
it here for completeness.
\begin{definition}
  \label{def:beef}
  The pair
  $(H,Q) \in \C2 (\reali^2;\reali) \times \C2 ([0,T]\times \overline
  \Omega \times \reali^2; \reali^N)$
  is called a \emph{boundary entropy--entropy flux pair} with respect
  to $f$ if
  \begin{enumerate}[i)]
  \item for all $w \in \reali$ the function $z \mapsto H (z,w)$ is
    convex;
  \item for all $t \in [0,T]$, $x\in\overline\Omega$ and
    $z,\, w \in \reali$,
    $\partial_3 Q (t,x,z,w) = \partial_1 H (z,w) \,\partial_3
    f\left(t,x,z\right)$;
  \item for all $t \in [0,T]$, $x\in\overline\Omega$ and
    $w \in \reali$, $H (w,w)= 0$, $Q (t,x,w,w) = 0$ and
    $\partial_1 H (w,w) = 0$.
  \end{enumerate}
\end{definition}

\noindent Note that if $H$ is as above, then $H \geq 0$.  \smallskip

We now extend the definition given by Otto
(see~\cite[Proposition~2]{Otto} and
also~\cite[Theorem~7.31]{MalekEtAlBook}) to account for non autonomous
fluxes and arbitrary source terms. The concept of boundary
entropy--entropy flux pairs introduced above characterises the
definition.
\begin{definition}
  \label{def:resol}
  A \emph{regular entropy solution (RE--solution)} to the
  initial--boundary value problem~\eqref{eq:1} on the interval $[0,T]$
  is a map $u \in \L\infty ([0,T] \times \Omega; \reali)$ such that
  for any boundary entropy--entropy flux pair $(H, Q)$, for any
  $k \in \reali$ and for any test function
  $\phi \in \Cc1 (]-\infty,T[ \times \reali^N; \reali_+)$
  \begin{equation}
    \label{eq:re}
    \begin{aligned}
      & \int_0^T \int_{\Omega} \left[ H \left(u (t,x), k
        \right) \partial_t \phi (t,x) + Q \left(t,x,u (t,x), k \right)
        \cdot \nabla \phi (t,x)\right] \d{x} \d{t}
      \\
      & + \int_0^T \int_{\Omega}
      \partial_1 H \left(u (t,x), k \right) \, \left[ F \left(t,x,u
          (t,x)\right) - \div f\left(t,x,u (t,x)\right) \right] \phi
      (t,x) \d{x} \d{t}
      \\
      & + \int_0^T \int_{\Omega} \div Q \left(t,x,u (t,x), k \right)
      \, \phi (t,x)\d{x} \d{t}
       + \int_{\Omega} H\left(u_o (x), k \right) \, \phi (0,x) \d{x}
      \\
      & + \norma{\partial_u f}_{\L\infty ([0,T] \times \Omega \times
        \mathcal{U}; \reali^N)} \int_0^T \int_{\partial \Omega} H \left(u_b
        (t,\xi), k \right) \, \phi (t,\xi) \d\xi \d{t} \geq 0,
    \end{aligned}
  \end{equation}
  where $\mathcal{U}$ is the interval
  $\mathcal{U}=[-U,U]$, with $U=\norma{u}_{\L\infty ([0,T] \times \Omega;\reali)}$.
\end{definition}

\noindent A comment on the constant appearing in the last line of the integral
inequality above is at the end of Section~\ref{sec:1d}.

\begin{remark}\label{rem:fnztest}
  {\rm Observe that an equivalent definition of solution can be
    obtained considering test functions
    $\phi \in \Cc1 (\reali \times \reali^N; \reali_+)$ and the
    following integral inequality:
 \begin{align*}
          & \int_0^T \int_{\Omega} \left[ H \left(u (t,x), k
        \right) \partial_t \phi (t,x) + Q \left(t,x,u (t,x), k \right)
        \cdot \nabla \phi (t,x)\right] \d{x} \d{t}
      \\
      & + \int_0^T \int_{\Omega}
      \partial_1 H \left(u (t,x), k \right) \, \left[ F \left(t,x,u
          (t,x)\right) - \div f\left(t,x,u (t,x)\right) \right] \phi
      (t,x) \d{x} \d{t}
      \\
      & + \int_0^T \int_{\Omega} \div Q \left(t,x,u (t,x), k \right)
      \, \phi (t,x)\d{x} \d{t}
      \\
      & + \int_{\Omega} H\left(u_o (x), k \right) \, \phi (0,x) \d{x}
      - \int_{\Omega}H\left(u (T,x), k \right) \, \phi (T,x) \d{x}
      \\
      & + \norma{\partial_u f}_{\L\infty ([0,T] \times \Omega \times
        \mathcal{U}; \reali^N)} \int_0^T \int_{\partial \Omega} H \left(u_b
        (t,\xi), k \right) \, \phi (t,\xi) \d\xi \d{t} \geq 0.
 \end{align*}  }
\end{remark}

A similar definition of solution is given by Vovelle
in~\cite[Definition~1]{Vovelle}, see also~\cite[Definition~1]{Martin},
using the so called \emph{Kru\v zkov semi-entropy--entropy flux pairs}, which are Lipschitz continuous functions, thus less regular than the
boundary entropies considered in the definition of
RE--solution.

\begin{definition}
  \label{def:mvsol}
  A \emph{semi-entropy solution (MV--solution)} to the initial--boundary value
  problem~\eqref{eq:1} on the interval $[0,T]$ is a map
  $u \in \L\infty ([0,T] \times \Omega; \reali)$ such that for any
  $k \in \reali$ and for any test function
  $\phi \in \Cc1 (]-\infty,T[ \times \reali^N; \reali_+)$
  \begin{align}
    \nonumber
    & \int_0^T \int_{\Omega} \left(u (t,x) - k \right)^\pm
      \, \partial_t \phi (t,x) \d{x} \d{t}
    \\
    \nonumber
    & + \int_0^T \int_{\Omega}  \sgn {}^\pm (u (t,x) -k) \;
      \left(f\left(t,x,u (t,x) \right) -f\left(t,x,k \right)\right)
      \cdot \nabla \phi (t,x) \d{x} \d{t}
    \\
    \label{eq:mv}
    & + \int_0^T \int_{\Omega} \sgn{}^\pm (u (t,x) -k)
      \left[
      F\left(t,x,u (t,x) \right) - \div f \left(t,x,k \right)
      \right] \, \phi (t,x) \d{x} \d{t}
    \\
    \nonumber
    & +
      \int_{\Omega} \left(u_o (x) - k \right)^\pm \, \phi (0,x) \d{x}
    \\
    \nonumber
    & +
      \norma{\partial_u f}_{\L\infty ([0,T] \times \Omega \times \mathcal{U}; \reali^N)}
      \int_0^T \int_{\partial \Omega}
      \left(u_b(t,\xi) - k\right)^\pm \, \phi (t,\xi) \d\xi \d{t} \geq 0,
  \end{align}
  where $\mathcal{U}$ is the interval $\mathcal{U}=[-U,U]$, with
  $U=\norma{u}_{\L\infty ([0,T] \times \Omega;\reali)}$.
\end{definition}

Before entering into the details of the link between these two
definitions of solution to~\eqref{eq:1}, we emphasise a feature they
share.
\begin{remark}
  \label{rem:stableL1}
  {\rm Both Definitions~\ref{def:resol} of RE--solution
    and~\ref{def:mvsol} of MV--solution are stable under
    $\L1$-convergence. This remarkable feature is underlined
    in~\cite[Chapter~2, Remark~7.33]{MalekEtAlBook} for the particular
    definition given by Otto, but it is immediate to see that it
    extends to both Definitions~\ref{def:resol}
    and~\ref{def:mvsol}. More precisely, let $u_o^n$ and $u_b^n$ be
    sequences of initial and boundary data converging in $\L1$ to
    $u_o$ and $u_b$ respectively. Let $u^n$ be a solution
    to~\eqref{eq:1}, according to either of the two definitions, with
    initial datum $u_o^n$ and boundary datum $u_b^n$. Assume that
    $u_n$ converges to $u$ in $\L1$. Then, this limit function $u$ is
    a solution to~\eqref{eq:1}, according to the same definition, with
    initial datum $u_o$ and boundary datum $u_b$.}
\end{remark}
\smallskip

Our first aim is to establish a connection between
Definition~\ref{def:resol} of RE--solution and
Definition~\ref{def:mvsol} of MV--solution. An intermediate step is
constituted by the following Lemma, which gives a link between the
boundary entropy--entropy flux pairs exploited in
Definition~\ref{def:resol} and the Kru\v zkov semi-entropy--entropy
flux pairs used in
Definition~\ref{def:mvsol}. 
\begin{lemma}[{\cite[Lemma~1]{Vovelle} and~\cite[Lemma~3]{Martin}}]
  \label{lem:vovelle}
  Let $\eta \in \C2 (\reali;\reali)$ be a convex function such that
  there exists $w \in [A,B]$ with $\eta(w) = 0$ and $\eta' (w)=0$.
  Then $\eta$ can be uniformly approximated on $[A,B]$ by applications
  of the kind
  \begin{displaymath}
    s \mapsto \sum_{i=1}^p \alpha_i \, (s-\kappa_i)^-
    +
    \sum_{j=1}^q \beta_j\, (s-\tilde \kappa_j)^+,
  \end{displaymath}
  where $\alpha_i \geq 0$, $\beta_j\geq 0$,
  $\kappa_i, \, \tilde \kappa_j \in [A,B]$.

  Conversely, there exists a sequence of boundary entropy--entropy
  flux pairs which converges to the Kru\v zkov semi-entropy--entropy
  flux pairs.
\end{lemma}
Thanks to Lemma~\ref{lem:vovelle}, the equivalence between
RE--solution and MV--solution follows immediately. For the detailed
proof we refer to Section~\ref{sec:TD}.
\begin{theorem}
  \label{thm:remv}
  Let $u \in \L\infty ([0,T] \times \Omega; \reali)$. Then $u$ is a
  RE--solution to~\eqref{eq:1}, in the sense of
  Definition~\ref{def:resol}, if and only if $u$ is a MV--solution
  to~\eqref{eq:1}, in the sense of Definition~\ref{def:mvsol}.
\end{theorem}

\section{Behaviour at the Boundary}
\label{sec:bc}

We now focus our attention on the way the boundary conditions are
fulfilled according to the definitions of solution to~\eqref{eq:1}
introduced in Section~\ref{sec:def}.  All the proofs are deferred to
Section~\ref{sec:TD}.

The following Lemma is a generalisation to problem~\eqref{eq:1}
of~\cite[Lemma~7.34]{MalekEtAlBook}. It states the way the boundary
conditions are satisfied in the case of a RE--solution
to~\eqref{eq:1}.

\begin{lemma}
  \label{lem:rebc}
  Let $u \in \L\infty ([0,T] \times \Omega; \reali)$ be a RE--solution
  to~\eqref{eq:1}, according to Definition~\ref{def:resol}. Then, for
  all boundary entropy--entropy flux pairs $(H,Q)$ and for all
  $\beta \in \L1 (]0,T[ \times \partial \Omega; \reali_+)$
  \begin{equation}
    \label{eq:rebc}
    \esslim_{\rho \to 0^+} \int_0^T \int_{\partial\Omega}
    Q \left(t, \xi, u\left(t, \xi - \rho \, \nu (\xi)\right), u_b (t,\xi)\right)
    \cdot \nu (\xi) \, \beta (t, \xi) \d\xi\d{t} \geq 0,
  \end{equation}
  $\nu (\xi)$ being the exterior normal to $\xi \in \partial \Omega$.
\end{lemma}

\begin{remark}
  {\rm Due to the equivalence between RE--solution and MV--solution
    proved in Theorem~\ref{thm:remv}, the boundary conditions are
    satisfied in the sense of~\eqref{eq:rebc} also in the case of a
    MV--solution to~\eqref{eq:1}. See~\cite[Lemma~4]{Martin}, and
    also~\cite[Remark~3]{Vovelle}, for a different proof
    of~\eqref{eq:rebc}, starting from Kru\v zkov semi-entropy--entropy
    flux pairs.}
\end{remark}

An alternative formulation of the boundary conditions is also
possible, both in the case of RE--solution and MV--solution,
see~\cite[Lemma~7.12]{MalekEtAlBook} and~\cite[Lemma~16]{Martin}.
\begin{lemma}
  \label{lem:mvbc}
  Let $u \in \L\infty ([0,T] \times \Omega; \reali)$ be a RE--solution
  (or MV--solution) to~\eqref{eq:1} in the sense of
  Definition~\ref{def:resol} (Definition~\ref{def:mvsol}).  Define the
  function
  $\mathcal{F} \in \C0 ([0,T] \times \overline \Omega \times \reali^3;
  \reali^N)$:
  \begin{equation}
    \label{eq:2}
    \mathcal{F} (t, x, z, w, k) =
    \begin{cases}
      f (t,x,w) - f (t,x,z) & \mbox{ for } z \leq w \leq k,
      \\
      0 & \mbox{ for } w \leq z \leq k,
      \\
      f (t,x,z) - f (t,x,k) & \mbox{ for } w \leq k \leq z,
      \\
      f (t,x,k) - f (t,x,z) & \mbox{ for } z \leq k \leq w,
      \\
      0 & \mbox{ for } k \leq z \leq w,
      \\
      f (t,x,z) - f (t,x,w) & \mbox{ for } k \leq w \leq z.
    \end{cases}
  \end{equation}
  Then, for all
  $\beta \in \L1 (]0,T[ \times \partial \Omega; \reali_+)$ and for all
  $k \in \reali$
  \begin{equation}
    \label{eq:mvbc}
    \esslim_{\rho \to 0^+} \int_0^T \int_{\partial\Omega}
    \mathcal{F} \left(t, \xi, u\left(t, \xi - \rho \, \nu (\xi)\right), u_b (t,\xi), k \right)
    \cdot \nu (\xi) \, \beta (t, \xi) \d\xi\d{t} \geq 0.
  \end{equation}
\end{lemma}

\begin{remark}
  {\rm Observe that the function $\mathcal{F}$ defined in~\eqref{eq:2} can
  be written also as follows
  \begin{align*}
    \mathcal{F} (t, x, z, w, k) = \
    & \dfrac12 \left[
      \sgn (z-w)\left(f (t,x,z) - f (t,x,w) \right) \right.
    \\[1pt]
    & \quad -
      \sgn (k-w)\left(f (t,x,k) - f (t,x,w) \right)
    \\[1pt]
    & \quad \left.
      +
      \sgn (z-k)\left(f (t,x,z) - f (t,x,k) \right)
      \right],
  \end{align*}
  and also
  \begin{align*}
    \mathcal{F} (t, x, z, w, k) = \
    & \sgnp\left(z-\max\{w,k\}\right)
      \left(f (t,x,z) - f\left(t,x,\max\{w,k\}\right)\right)
    \\
    & + \sgnm\left(z-\min\{w,k\}\right)
      \left(f (t,x,z) - f\left(t,x,\min\{w,k\}\right)\right).
  \end{align*}}
\end{remark}

\section{Solutions with Traces}
\label{sec:solTR}

So far we have considered two definitions of solution to~\eqref{eq:1},
sought in $\L\infty$, and proved their equivalence. The RE--definition
involves regular entropies, while the MV--definition deals with
Lipschitz continuous ones. In this Section we present two additional
definitions of solution to~\eqref{eq:1}, in which the trace of the
solution at the boundary appears explicitly. The idea is to draw a
parallel with RE and MV--solutions: indeed, the two definitions we are
going to introduce are characterised by regular and Lipschitz
continuous entropies respectively, and we prove that they are
equivalent.

\smallskip

Since the existence of the trace of the solution is required, more
regularity is needed on the solution with respect to
Definitions~\ref{def:resol} of RE--solutions and~\ref{def:mvsol} of
MV--solutions. To this aim, introduce the following space:
\begin{definition}
  \label{def:trspace}
  A function $u$ belongs to the space
  $\mathcal{TR}^\infty ([0,T] \times \Omega;\reali)$ if there exists a
  function $ \tr u \in \L\infty ([0,T] \times \partial\Omega; \reali)$
  such that
  \begin{equation}
    \label{eq:trspace}
    \esslim_{r \to 0^+} \int_0^T \int_{\partial\Omega}
    \modulo{u \left(t,\xi -r\, \nu (\xi)\right) - \tr u (t, \xi)} \d\xi \d{t}
    =0.
  \end{equation}
\end{definition}

We remark the following. Bardos, le Roux and N\'ed\'elec
in~\cite{BardosLerouxNedelec} consider solutions in
$\BV ([0,T] \times \Omega;\reali) \subset \mathcal{TR}^\infty ([0,T]
\times \Omega;\reali)$:
indeed, $\BV$ functions admit a trace at the boundary reached by $\L1$
convergence, see~\cite[Lemma~1]{BardosLerouxNedelec},
\cite[Paragraph~5.3]{EvansGariepy}, \cite[Chapter~2]{GiustiBook}
and~\cite[Appendix]{bordo}. Consider now the case of $\L\infty$
solutions, which are functions $u$ satisfying in the sense of
distribution on $]0,T[ \times \pint \Omega$ the inequality
\begin{displaymath}
  \partial_t \eta (u) + \div q (t,x,u) \leq 0,
\end{displaymath}
for any $(\eta,q)$ (classical) entropy--entropy flux pair (with
respect to $f$, see Definition~\ref{def:eef}). Panov proves
in~\cite[Theorem~1.4 and Remark~6.3]{Panov2007} the existence of the
trace at the boundary for $\L\infty$ solutions to~\eqref{eq:1}, under
the following non-degeneracy condition on the flux: the function $f$
is continuous and such that for
a.e.~$(t,\xi) \in \reali_+ \times \partial\Omega$ and all
$(s,y) \in \reali_+ \times \reali^N \setminus \ \{(0,0)\}$ the
function $u \to s \, u + y \, f (t,\xi,u)$ is not constant on
non-degenerate intervals, i.e.~$f$ satisfies the following genuine non
linearity condition
\begin{displaymath}
  \mathcal{L} \left(\left\{
      u | s + y \, \partial_u f (t,\xi,u) = 0
    \right\}\right) = 0
  \quad \mbox{ for every }
  (s,y) \in \reali_+ \times \reali^N \setminus \ \{(0,0)\},
\end{displaymath}
where $\mathcal{L}$ is the Lebesgue measure and
$(t,\xi) \in \reali_+ \times \partial\Omega$. As pointed out also
in~\cite{Vasseur}, the above assumption allows to avoid flux functions
whose restriction to an open subset is linear.


\smallskip

The following definition uses the (classical) \emph{entropy--entropy
  flux pairs}, see Definition~\ref{def:eef}. It extends the particular
case of scalar conservation laws with autonomous fluxes considered
in~\cite[Chapter~6, Definition~6.9.1]{DafermosBook} and
in~\cite[Chapter~15]{Serre2}, see also~\cite[Definition~2.5]{bordo}.
\begin{definition}
  \label{def:esol}
  An \emph{entropy solution (E--solution)} to the initial--boundary
  value problem~\eqref{eq:1} on the interval $[0,T]$ is a map
  $u \in (\L\infty \cap \mathcal{TR}^\infty )([0,T] \times \Omega;
  \reali)$
  such that for any entropy--entropy flux pair $(\eta, q)$ and for any
  test function $\phi \in \Cc1 (]-\infty,T[ \times \reali^N; \reali_+)$
  \begin{align}
    \nonumber
    & \int_0^T \int_{\Omega} \left[ \eta\left(u (t,x)\right)
      \, \partial_t \phi (t,x)
      + q\left(t,x,u (t,x) \right) \cdot \nabla \phi (t,x)\right] \d{x} \d{t}
    \\
    \nonumber
    & + \int_0^T \int_{\Omega} \eta'\left(u (t,x)\right)
      \left[
      F\left(t,x,u (t,x) \right) - \div f \left(t,x,k \right)
      \right] \, \phi (t,x) \d{x} \d{t}
    \\
    \label{eq:e}
    & + \int_0^T \int_{\Omega} \div q \left(t,x,u (t,x) \right)
      \, \phi (t,x) \d{x} \d{t}
    \\
    \nonumber
    & +
      \int_{\Omega} \!\eta\left(u_o (x)\right) \phi (0,x) \d{x}
      -
      \int_0^T\! \int_{\partial \Omega}\!
      q \left(t,\xi,u_b (t,\xi) \right) \cdot
      \nu (\xi) \, \phi (t,\xi) \d\xi \d{t}
    \\
    \nonumber
    & +
      \int_0^T \int_{\partial \Omega}
      \eta'\left(u_b(t,\xi)\right)
      \left(
      f\left(t,\xi, u_b(t,\xi)\right) - f (t,\xi,\tr u (t,\xi))\right)
      \cdot
      \nu (\xi) \, \phi (t,\xi) \d\xi \d{t} \geq 0.
  \end{align}
\end{definition}

We now recall the definition of solution to~\eqref{eq:1} due to
Bardos, le Roux and N\'ed\'elec~\cite[p.~1028]{BardosLerouxNedelec},
which exploits the classical \emph{Kru\v zkov entropy--entropy flux
  pairs}.
\begin{definition}
  \label{def:blnsol}
  A \emph{Kru\v zkov-entropy solution (BLN--solution)} to the
  initial--boundary value problem~\eqref{eq:1} on the interval $[0,T]$
  is a map
  $u \in (\L\infty \cap \mathcal{TR}^\infty )([0,T] \times \Omega;
  \reali)$
  such that for any $k \in \reali$ and for any test function
  $\phi \in \Cc1 (]-\infty,T[ \times \reali^N; \reali_+)$
  \begin{align}
    \nonumber
    & \int_0^T \int_{\Omega} \modulo{u (t,x) - k}
      \, \partial_t \phi (t,x) \d{x} \d{t}
    \\
    \nonumber
    & + \int_0^T \int_{\Omega}  \sgn (u (t,x) -k) \;
      \left(f\left(t,x,u (t,x) \right) -f\left(t,x,k \right)\right)
      \cdot \nabla \phi (t,x) \d{x} \d{t}
    \\
    \label{eq:bln}
    & + \int_0^T \int_{\Omega} \sgn (u (t,x) -k)
      \left[
      F\left(t,x,u (t,x) \right) - \div f \left(t,x,k \right)
      \right] \, \phi (t,x) \d{x} \d{t}
    \\
    \nonumber
    & +
      \int_{\Omega} \modulo{u_o (x) - k }\, \phi (0,x) \d{x}
    \\
    \nonumber
    & -
      \int_0^T \int_{\partial \Omega}
      \sgn\left(u_b (t,\xi) -k\right)
      \left(f\left(t,\xi,\tr u(t,\xi)\right) - f (t,\xi,k)\right) \cdot
      \nu (\xi) \, \phi (t,\xi) \d\xi \d{t} \geq 0.
  \end{align}
\end{definition}

E--solutions, as in Definition~\ref{def:esol}, and BLN--solutions, as
in Definition~\ref{def:blnsol}, are actually equivalent,
see~\cite[Proposition~2.6]{bordo}. The proof of the equivalence
between these two Definitions of solution is based on an analogous of
Lemma~\ref{lem:vovelle} and is briefly sketched in
Section~\ref{sec:TD}.
\begin{theorem}
  \label{thm:ebln}
  The map
  $u \in (\L\infty \cap \mathcal{TR}^\infty )([0,T] \times \Omega;
  \reali)$
  is an E--solution to~\eqref{eq:1}, in the sense of
  Definition~\ref{def:esol}, if and only if $u$ is a BLN--solution
  to~\eqref{eq:1}, in the sense of Definition~\ref{def:blnsol}.
\end{theorem}

Before studying the relation among all the considered definitions, we
provide the analogous to Lemma~\ref{lem:rebc} and
Lemma~\ref{lem:mvbc}, explaining the way E--solutions and
BLN--solutions to~\eqref{eq:1} fulfil the boundary
conditions. Concerning E--solutions, the following Lemma holds.
\begin{lemma}
  \label{lem:ebc}
  Let
  $u \in (\L\infty \cap \mathcal{TR}^\infty)([0,T] \times \Omega;
  \reali)$
  be an E--solution to~\eqref{eq:1} in the sense of
  Definition~\ref{def:esol}. Then, for all (classical)
  entropy--entropy flux pairs $(\eta,q)$ and for
  a.e.~$(t,\xi) \in ]0,T[ \times \partial \Omega$,
  \begin{equation}
    \label{eq:ebc}
    \left[
      q\left(t,\xi,\tr u (t,\xi)\right)
      \!-
      q\left(t,\xi, u_b (t,\xi)\right)
      \!- \eta'\left(u_b (t,\xi)\right)\!\!
      \left(
        \!f\!\left(t,\xi,\tr u (t,\xi)\right)
        \!-\!
        f \!\left(t,\xi, u_b (t,\xi)\right)\!
      \right)\right] \cdot \nu (\xi) \geq 0.
  \end{equation}
\end{lemma}
The proof is deferred to Section~\ref{sec:TD}.
Observe that condition~\eqref{eq:ebc} is the generalisation to the
multidimensional case of the \emph{boundary entropy inequality} due to
Dubois and LeFloch~\cite[Theorem~1.1]{DuboisLeFloch}.

In the following Lemma we recall the well-known BLN condition, linking the
boundary datum and the trace of the solution. The proof is in
Section~\ref{sec:TD}.
\begin{lemma}
  \label{lem:blnbc}
  Let
  $u \in (\L\infty \cap \mathcal{TR}^\infty)([0,T] \times \Omega;
  \reali)$
  be a BLN--solution to~\eqref{eq:1} in the sense of
  Definition~\ref{def:blnsol}. Then, for all $k \in \reali$ and for
  a.e.~$(t,\xi) \in ]0,T[ \times \partial \Omega$,
  \begin{equation}
    \label{eq:6}
    \left(\sgn\left(\tr u (t,\xi) - k\right)
      -
      \sgn\left(u_b (t,\xi) - k\right)
    \right)
    \left(
      f\left(t,\xi,\tr u (t,\xi)\right)
      -
      f (t,\xi,k)
    \right) \cdot \nu (\xi) \geq 0.
  \end{equation}
  Moreover, condition~\eqref{eq:6} is equivalent to the following: for
  all $k\in \mathcal{I}[\tr u (t,\xi), u_b (t,\xi)]$ and
  a.e.~$(t,\xi) \in ]0,T[ \times \partial \Omega$
  \begin{equation}
    \label{eq:blnbc}
    \sgn\left(\tr u (t,\xi) - u_b (t,\xi)\right)
    \left(
      f\left(t,\xi,\tr u (t,\xi)\right)
      -
      f (t,\xi,k)
    \right) \cdot \nu (\xi) \geq 0.
  \end{equation}
\end{lemma}

The following Proposition constitutes the basis for the proof of the
equivalence of all the definitions of solution to~\eqref{eq:1}
presented so far. It is a generalisation
of~\cite[Lemma~7.24]{MalekEtAlBook} to problem~\eqref{eq:1}: it takes
into account non autonomous fluxes and arbitrary source terms. In
particular, this Proposition provides a connection among the ways the
boundary conditions are understood according to the various
definitions of solution introduced so far. However, we need to require
the existence of the trace of the solution at the boundary. For
further details about the trace, see the references at the beginning
of Section~\ref{sec:solTR}. The proof is deferred to
Section~\ref{sec:TD}.

\begin{proposition}
  \label{prop:equivbc}
  Let $u_b \in \L\infty ([0,T] \times \partial \Omega;\reali)$ and
  $u \in (\L\infty \cap \mathcal{TR}^\infty) ([0,T] \times
  \Omega;\reali)$. Then the following statements are equivalent:
  \begin{enumerate}
  \item\label{item:1} ~\eqref{eq:rebc} holds, for any boundary
    entropy--entropy flux pair $(H,Q)$ and for any
    $\beta \in \L1 (]0,T[ \times \partial \Omega;\reali_+)$;

  \item\label{item:2} ~\eqref{eq:mvbc} holds, for any
    $\beta \in \L1 (]0,T[ \times \partial \Omega;\reali_+)$ and for
    all $k\in \reali$;

  \item\label{item:4} for
    a.e.~$(t,\xi) \in ]0,T[ \times \partial \Omega$ and for all
    $k \in \reali$ it holds
    \begin{equation}
      \label{eq:5}
      \mathcal{F}\left(t, \xi, \tr u (t, \xi), u_b (t,\xi), k\right)
      \cdot \nu (\xi)\geq 0,
    \end{equation}
    with $\mathcal{F}$ as in~\eqref{eq:2};

  \item\label{item:5} ~\eqref{eq:blnbc} holds for
    a.e.~$(t,\xi) \in ]0,T[ \times \partial \Omega$ and for all
    $k \in \mathcal{I}[\tr u (t,\xi), u_b (t,\xi)]$;

  \item\label{item:6} ~\eqref{eq:ebc} holds for
    a.e.~$(t,\xi) \in ]0,T[ \times \partial \Omega$ and for any
    entropy--entropy flux pair $(\eta,q)$;

  \item\label{item:3} for
    a.e.~$(t,\xi) \in ]0,T[ \times \partial \Omega$ and all
    entropy--entropy flux pair $(\eta, q)$, such that
    $\eta' \left(u_b (t,\xi)\right)\! = 0$ and
    $q\left(t, \xi, u_b (t,\xi)\right)=0$, it holds
    \begin{equation}
      \label{eq:dubois}
      q\left(t,\xi, \tr u (t,\xi)\right) \cdot \nu (\xi) \geq 0.
    \end{equation}
  \end{enumerate}
\end{proposition}



We can now state our main result: given that $u$ admits a trace in the
sense of~\eqref{eq:trspace}, we prove that the Definitions of solution
presented in this Section, that is E--solution and BLN--solution, are
equivalent to the Definitions of solution introduced in
Section~\ref{sec:def}, i.e.~RE--solution and MV--solution.
\begin{theorem}
  \label{thm:mvbln}
  Let
  $u \in (\L\infty \cap \mathcal{TR}^\infty)([0,T] \times \Omega;
  \reali)$.
  Then $u$ is a RE--solution to~\eqref{eq:1} according to
  Definition~\ref{def:resol}, or equivalently a MV--solution
  to~\eqref{eq:1} in the sense of Definition~\ref{def:mvsol}, if and
  only if $u$ is an E--solution to~\eqref{eq:1} according to
  Definition~\ref{def:esol}, or equivalently a BLN--solution
  to~\eqref{eq:1} in the sense of Definition~\ref{def:blnsol}.
\end{theorem}
\noindent The proof is deferred to Section~\ref{sec:TD}. Remark that,
according to the results by Panov~\cite{Panov2007} recalled at the
beginning of this section, $\L\infty$ solutions admit a trace at the
boundary in the case of non-degenerate fluxes, thus in those cases
there is no need to consider the intersection with the space
$\mathcal{TR}^\infty([0,T] \times \Omega; \reali)$.

\section{Strong Solutions}
\label{sec:strongsol}

For completeness, we recall below the definition of strong (smooth)
solution to~\eqref{eq:1}.

\begin{definition}
  \label{def:strong}
  A \emph{strong solution} to the initial--boundary value
  problem~\eqref{eq:1} on the interval $[0,T]$ is a map
  $u \in \C1 ([0,T]\times \Omega;\reali) \cap \C0([0,T] \times
  \overline \Omega;\reali)$
  which satisfies pointwise the equation and the initial condition,
  and it is such that, for all
  $(t,\xi)\in]0,T[ \times \partial \Omega$ and for all
  $k \in \mathcal{I}[u (t,\xi), u_b (t,\xi)]$,
  \begin{equation}
    \label{eq:strongbc}
    \sgn\left(u (t,\xi) - u_b (t,\xi)\right)
    \left(
      f\left(t,\xi,u (t,\xi)\right)
      -
      f (t,\xi,k)
    \right) \cdot \nu (\xi) \geq 0.
  \end{equation}
\end{definition}

Note that condition~\eqref{eq:strongbc} reduces to~\eqref{eq:blnbc}:
the difference is that strong solutions are defined up to the
boundary, and therefore the notion of trace is not needed. For further
details on the boundary conditions for smooth solution, including an
heuristic derivation, see~\cite[Chapter~2, Section~6]{MalekEtAlBook}.

The following result holds.
\begin{proposition}
  \label{propo:strongre}
  Let
  $u \in \C1 ([0,T]\times \Omega;\reali) \cap \C0([0,T] \times
  \overline \Omega;\reali)$
  be a strong solution to~\eqref{eq:1} in the sense of
  Definition~\ref{def:strong}. Then $u$ is also a $RE$--solution
  to~\eqref{eq:1}, in the sense of Definition~\ref{def:mvsol}.
\end{proposition}
\noindent Obviously, due to the equivalence among the definitions of
solution proven in Theorem~\ref{thm:mvbln}, every strong solution is
also a MV--solution, an E--solution and a BLN--solution.

The proof follows the line of the second part of the proof of
Theorem~\ref{thm:mvbln} and it is hence omitted. The main difference
is that the solution itself at the boundary is considered, instead of
its trace.

\section{The 1 Dimensional Case}
\label{sec:1d}

In this section we focus on the case $N=1$, i.e.~$\Omega$ is the
interval $]a,b[$, with $a, b \in \reali$. The boundary datum is assigned at the end points of the interval: for
$t \in \reali_+$
\begin{align*}
  u(t,a) = \
  & u_b (t,a),
  & u (t,b)= \
  & u_b (t,b).
\end{align*}
We write explicitly how the last line of the integral inequality in
the definition of solution reads in the case of the RE--definition and
of the E--definition, the other two cases being completely
analogous. Observe that the exterior normal to $\partial \Omega$ in
$a$ is $-1$, while in $b$ is $+1$.
\begin{itemize}
\item {\bf RE--definition:}
  \begin{displaymath}
    + \norma{\partial_u f}_{\L\infty ([0,T] \times \Omega \times \mathcal{U};\reali)}
    \int_0^T\left[
      H \left(u_b(t,a), k \right) \, \phi (t,a)
      +
      H \left(u_b(t,b), k \right) \, \phi (t,b)
    \right]\d{t}.
  \end{displaymath}

\item {\bf E--definition:}
  \begin{multline*}
    + \int_0^T q \left(t,a,u_b (t,a)\right) \, \phi (t,a)\d{t}
    -
    \int_0^T q \left(t,b,u_b (t,b)\right) \, \phi (t,b)\d{t}
    \\
    -
    \int_0^T \eta'\left(u_b (t,a)\right)
    \left(f\left(t,a,u_b (t,a)\right) - f\left(t,a,\tr u (t,a)\right)\right)
    \phi (t,a) \d{t}
    \\
    +
    \int_0^T \eta'\left(u_b (t,b)\right)
    \left(f\left(t,b,u_b (t,b)\right) - f\left(t,b,\tr u (t,b)\right)\right)
    \phi (t,b) \d{t}.
  \end{multline*}
\end{itemize}
What is immediately evident is the presence of the minus sign in the
last case, while the first contains only sums. The minus sign is due
to the scalar product with the exterior normal to $\partial \Omega$,
which occurs in the E--definition, and in the BLN--definition as
well. It can be seen that there is no need for a minus sign neither in
the RE--definition nor in the MV--definition.

\medskip

We can exploit the one dimensional setting to analyse a feature of the
RE--definition and the MV--definition. Indeed, the integral
inequalities of these two definitions involve the constant
$\norma{\partial_u f}_{\L\infty ([0,T] \times \Omega \times
  \mathcal{U};\reali^N)}$,
where $\mathcal{U}=[-U,U]$, with
$U= \norma{u}_{\L\infty ([0,T] \times \Omega;\reali)}$. This is
nothing but the Lipschitz constant of $f$ with respect to $u$,
therefore one may wonder whether it is possible to consider a
different constant, either smaller or larger, and to still get a
solution. One can see, through the following one dimensional example,
that the above constant is indeed the smallest possible.  Consider the
following problem
\begin{equation*}
  \left\{
    \begin{array}{l@{\qquad}r@{\,}c@{\,}l}
      \partial_t u (t, x)
      +
      \partial_x f \left(u (t,x)\right)
      =
      0
      & (t,x)
      & \in
      & \reali_+ \times ]a,b[
      \\
      u (0, x) = 1
      & x
      & \in
      & ]a,b[
      \\
      u (t,a) = u_b (t,a) = 1
      & t
      & \in
      & \reali_+
      \\
      u (t,b) = u_b (t,b) = -1
      & t
      & \in
      & \reali_+
    \end{array}
  \right.
  \qquad
  \mbox{ with }
  f (u) = \frac{u^2}{2}.
\end{equation*}
The solution is constant and equal to $1$. Fix $T>0$. Observe that
$\norma{u}_{\L\infty ([0,T] \times [a,b];\reali)}=1$. Consider the
integral inequality~\eqref{eq:mv} of the MV--definition, and use the
positive constant $c$ instead of
$\norma{\partial_u f}_{\L\infty ([0,T] \times [a,b] \times
  \mathcal{U};\reali)} = 1$.
Since we know that $u =1$ is the solution, simple computations show
that $c$ should be grater or equal to $1$.

Obviously, given a solution $u$, choosing a greater value of the
Lipschitz constant still ensures that $u$ is a solution: indeed the
constant is multiplied by a non negative term, both in the
RE and in the MV--definition.

\section{Notes on the Existence of Solutions}
\label{sec:existence}

We recap the results present in the literature
concerning the existence of solutions to problem~\eqref{eq:1},
specifying in each case the considered definition of solution.

As far as it concerns the case of an autonomous scalar conservation
laws, i.e.~$f=f (u)$ and $F=0$, Otto proves the existence and
uniqueness of a RE--solution to~\eqref{eq:1}, see~\cite{Otto} and
also~\cite[Chapter~2]{MalekEtAlBook} for more detailed proofs.

In~\cite{Martin}, Martin proves the existence and uniqueness of a
MV--solution to the general problem~\eqref{eq:1}, though imposing on the flux
and the source terms a condition which leads to a (simple) maximum
principle, see~\cite[Assumption~1.(iv)]{Martin}.

Bardos, le Roux and N\'ed\'elec prove in~\cite{BardosLerouxNedelec}
existence and uniqueness of a BLN--solution to~\eqref{eq:1} in the
case of homogeneous boundary conditions, i.e.~$u_b=0$. The result is
then generalised in~\cite{bordo} to allow for (regular) non zero
boundary data. The case of an autonomous scalar conservation laws with
BLN--definition is considered also by Serre
in~\cite[Chapter~15]{Serre2}. Dafermos
in~\cite[Chapter~6]{DafermosBook} focuses on the same autonomous
problem, although under homogeneous boundary conditions, and exploits
the RE--definition of solution.

\section{Technical Details}
\label{sec:TD}

\begin{proofof}{Theorem~\ref{thm:remv}}\vspace{-4mm}
  \paragraph{A RE--solution is a MV--solution.}
  It is enough to consider the following sequences of boundary
  entropy--entropy flux pairs
  \begin{align*}
    H_n (z,k) = \
    &
      \left(\left((z- k)^\pm\right)^2 + \frac{1}{n^2}\right)^{1/2} -\frac1n
    \\
    Q_n (t,x,z,k) = \
    &
      \int_k^z \partial_1 H_n (w, k) \, \partial_u f (t,x,w)\d{w}.
  \end{align*}
  Indeed, as $n$ goes to $+\infty$,
  \begin{align*}
    H_n (z,k) \longrightarrow \ & (z-k)^\pm
    \\
    Q_n (t,x,z,k) \longrightarrow \ & \sgn(z-k)^\pm \, \left(f (t,x,z) - f (t,x,k)\right),
  \end{align*}
  so that, in the limit,~\eqref{eq:re} becomes
  \begin{align*}
    & \int_0^T \int_{\Omega}
      \left(u (t,x) - k\right)^\pm \partial_t \phi (t,x) \d{x} \d{t}
      \\
    & + \int_0^T \int_{\Omega}
       \sgn\left( u (t,x)-k\right)^\pm \,
      \left(f \left(t,x,u (t,x)\right) - f (t,x,k)\right)
      \cdot \nabla \phi (t,x) \d{x} \d{t}
    \\
    & + \int_0^T \int_{\Omega}
      \sgn\left( u (t,x)-k\right)^\pm \,
      \left( F \left(t,x,u (t,x)\right)
      - \div f\left(t,x,u (t,x)\right)   \right) \phi (t,x) \d{x} \d{t}
    \\
    & + \int_0^T \int_{\Omega}
      \sgn\left( u (t,x)-k\right)^\pm
      \left(\div f \left(t,x,u (t,x) \right) -
      \div f \left(t,x, k \right)\right)
      \, \phi (t,x)\d{x} \d{t}
    \\
    & + \int_{\Omega} \left(u_o (x)- k \right)^\pm \, \phi (0,x) \d{x}
    \\
    & + \norma{\partial_u f}_{\L\infty ([0,T] \times \Omega \times \mathcal{U}; \reali^N)}
      \int_0^T \int_{\partial \Omega}  \left(u_b (t,\xi)-  k \right)^\pm \, \phi
      (t,\xi) \d\xi \d{t} \geq 0,
  \end{align*}
  where $\mathcal{U}=[-U,U]$, with
  $U=\norma{u}_{\L\infty ([0,T]\times \Omega;\reali)}$.  Combining the
  second and the third lines in the inequality above yields
  exactly~\eqref{eq:mv}.

  \paragraph{A MV--solution is a RE--solution.} Since
  $u \in \L\infty ([0,T]\times \Omega;\reali)$, there exist
  $A, \, B \in \reali$, with $A < B$, such that
  $A \leq u (t,x) \leq B$ for a.e.~$(t,x) \in [0,T]\times \Omega$, and
  hence $u \in \L\infty ([0,T]\times \Omega;[A,B])$. We can then apply
  Lemma~\ref{lem:vovelle}: each boundary entropy--entropy flux pair is
  uniformly approximated by a linear combination with positive
  coefficients of Kru\v zkov semi-entropy--entropy flux pairs. Thus
  the inequality in~\eqref{eq:mv} is preserved and~\eqref{eq:re} holds.
\end{proofof}

\begin{proofof}{Lemma~\ref{lem:rebc}}
  The proof extends~\cite[Lemma~7.34]{MalekEtAlBook} to consider non
  autonomous fluxes and general source terms.

  Let $(H,Q)$ be a boundary entropy--entropy flux pair,
  $k \in \reali$. The analogous of~\cite[Lemma~7.34,
  hypothesis~(7.35)]{MalekEtAlBook} is the following:
  \begin{align}
    \nonumber
    & \int_0^T \int_{\Omega} \left[ H \left(u (t,x), k
      \right) \partial_t \phi (t,x) + Q \left(t,x,u (t,x), k \right)
      \cdot \nabla \phi (t,x)\right] \d{x} \d{t}
    \\
    \label{eq:7.35}
    & + \int_0^T \int_{\Omega}
      \partial_1 H \left(u (t,x), k \right) \, \left[ F \left(t,x,u
      (t,x)\right) - \div f\left(t,x,u (t,x)\right) \right] \phi
      (t,x) \d{x} \d{t}
    \\
    \nonumber
    & + \int_0^T \int_{\Omega} \div Q \left(t,x,u (t,x), k \right)
      \, \phi (t,x)\d{x} \d{t}
      \geq 0,
  \end{align}
  which follows directly from the Definition~\ref{def:resol} of
  RE--solution when considering a test function
  $\phi \in \Cc1 (]0,T[ \times \pint \Omega; \reali_+)$.  Similarly,
  the analogous of~\cite[Lemma~7.34, hypothesis~(7.36)]{MalekEtAlBook}
  is the following:
  \begin{align}
    \nonumber
    & \int_0^T \int_{\Omega} \left[ H \left(u (t,x), k
      \right) \partial_t \phi (t,x) + Q \left(t,x,u (t,x), k \right)
      \cdot \nabla \phi (t,x)\right] \d{x} \d{t}
    \\
    \label{eq:7.36}
    & + \int_0^T \int_{\Omega}
      \partial_1 H \left(u (t,x), k \right) \, \left[ F \left(t,x,u
      (t,x)\right) - \div f\left(t,x,u (t,x)\right) \right] \phi
      (t,x) \d{x} \d{t}
    \\
    \nonumber
    & + \int_0^T \int_{\Omega} \div Q \left(t,x,u (t,x), k \right)
      \, \phi (t,x)\d{x} \d{t}
    \\ \nonumber
    & + \norma{\partial_u f}_{\L\infty ([0,T] \times \Omega \times \mathcal{U}; \reali^N)}
      \int_0^T \int_{\partial \Omega} H\left(u_b (t,\xi), k \right) \,
      \phi (t,\xi) \d\xi \d{t}
      \geq 0,
  \end{align}
  where $\mathcal{U}=[-U,U]$ with
  $U=\norma{u}_{\L\infty ([0,T] \times \Omega;\reali)}$, which follows
  directly from the Definition~\ref{def:resol} of RE--solution when
  considering a test function
  $\phi \in \Cc1 (]0,T[ \times \Omega; \reali_+)$.

  We proceed as in the proof of~\cite[Lemma~7.34]{MalekEtAlBook}. For
  the sake of simplicity, we restrict ourselves to the case of a
  half--space, i.e.
  \begin{align*}
    \Omega = \
    & \left\{
      x = (x', s) \in \reali^{N-1}\times \reali\colon  y<0
      \right\},
    \\
    \boldsymbol\nu = \
    & (0, \ldots, 0, 1) \in \reali^N,
    \\
    \Gamma = \
    & ]0,T[ \times \reali^{N-1}, \qquad r = \ (t,x') \in \Gamma,
    \\
    Q_T = \
    & \left\{
      p = (r,s) \colon r \in \Gamma , \, s<0
      \right\}.
  \end{align*}
  The general case can then be obtained by a covering argument, i.e.~by
  considering that the boundary $\partial \Omega$ can be locally
  replaced by the border of a half-space.

  For any boundary entropy--entropy flux pair $(H,Q)$, denote, for
  $(t,x) \in [0,T] \times \Omega$, $w \in \mathbb{Q}$,
  \begin{displaymath}
    \eta (z) = H (z,w)\,,
    \qquad
    q (t,x,z) = Q (t,x,z,w)\,.
  \end{displaymath}
  Choosing $\phi (t,x) = \phi (r,s) = \beta (r)\,\alpha (s)$, with
  $\alpha \in \Cc1 (]-\infty,0[;\reali_+)$, in~\eqref{eq:7.35} yields
  \begin{equation}
    \label{eq:7.37}
    -\int_{-\infty}^0 \int_\Gamma
    q \left( r, s, u (r,s) \right) \cdot \boldsymbol \nu \,
    \beta (r) \d{r} \, \alpha' (s)\d{s}
    \leq
    C \int_{-\infty}^0 \alpha (s) \d{s},
  \end{equation}
  where
  \begin{align*}
    C = \
    &
      \norma{\eta (u)}_{\L\infty (Q_T;\reali)}
      \int_{\Gamma} \modulo{\partial_t \beta (r)}\d{r}
      +
      \norma{q(\cdot,\cdot,u)}_{\L\infty (Q_T;\reali^N)}
      \int_\Gamma \modulo{\nabla_{x'} \beta(r)}\d{r}
    \\
    & +
      \left[
      \norma{\eta' (u)}_{\L\infty (Q_T;\reali)}
      \norma{(F - \div f)(\cdot, \cdot, u)}_{\L\infty (Q_T;\reali)}
      +
      \norma{\div q (\cdot, \cdot, u)}_{\L\infty (Q_T;\reali)}
      \right]
      \!\int_\Gamma \modulo{\beta (r)}\d{r}.
  \end{align*}
  Thanks to integration by parts on the left hand side
  of~\eqref{eq:7.37} and to the fact that $\alpha \geq 0$, we obtain
  that the function
  \begin{equation}
    \label{eq:7.38}
    s \mapsto \int_\Gamma
    q \left( r, s, u (r,s) \right) \cdot \boldsymbol \nu \,
    \beta (r) \d{r} - C \, s
  \end{equation}
  is non increasing on $]-\infty,0[$. Moreover, we have
  \begin{equation}
    \label{eq:4}
    \esslim_{s \to 0^-} \inf
    \int_\Gamma
    q \left( r, s, u (r,s) \right) \cdot \boldsymbol \nu \,
    \beta (r) \d{r}
    \geq
    -
    \esssup_{Q_T} \modulo{q (\cdot,\cdot,u)} \int_\Gamma \beta (r)\d{r}.
  \end{equation}
  Monotonicity~\eqref{eq:7.38} and boundedness from below~\eqref{eq:4}
  imply that the following quantity is finite
  \begin{equation}
    \label{eq:7.39}
    \esslim_{s \to 0^-}
    \int_\Gamma
    q \left( r, s, u (r,s) \right) \cdot \boldsymbol \nu \,
    \beta (r) \d{r}.
  \end{equation}
  From~\eqref{eq:7.36}, for all $\alpha\in\Cc1 (\reali;\reali_+)$ we
  get
  \begin{align*}
    &
      -\int_{-\infty}^0 \int_\Gamma
      q \left( r, s, u (r,s) \right) \cdot \boldsymbol \nu \,
      \beta (r) \d{r} \, \alpha' (s)\d{s}
    \\
    \leq \
    &  C \int_{-\infty}^0 \alpha (s) \d{s}
      +
      \norma{\partial_u f}_{\L\infty ([0,T] \times \Omega \times \mathcal{U}; \reali^N)}
      \int_\Gamma
      \eta \left(u_b (r)\right) \, \beta (r) \d{r} \alpha (0).
  \end{align*}
  Choose $\alpha_n (s) = (n \, s +1) \caratt{]-1/n,0[}$, mollify it
  properly and insert it in~\eqref{eq:7.39}: in the limit
  $n \to \infty$ we obtain
  \begin{equation}
    \label{eq:utile}
    \esslim_{s \to 0^-}
    \int_\Gamma
    q \left( r, s, u (r,s) \right) \cdot \boldsymbol \nu \,
    \beta (r) \d{r}
    \geq
    - \norma{\partial_u f}_{\L\infty ([0,T] \times \Omega \times \mathcal{U}; \reali^N)}
    \int_\Gamma
    \eta \left(u_b (r)\right) \, \beta (r) \d{r}.
  \end{equation}
  Let $J \subseteq \Cc1 (\Gamma;\reali_+)$ be a countable set of
  functions such that for all $\beta \in \L1 (\Gamma;\reali_+)$ there
  is a sequence $(\beta_n)$ in $J$ such that
  $\lim_{n} \beta_n = \beta$ in $\L1 (\Gamma;\reali_+)$. Therefore,
  \begin{displaymath}
    \lim_n \int_\Gamma
    q \left( r, s, u (r,s) \right) \cdot \boldsymbol \nu \,
    \beta_n (r) \d{r}
    =
    \int_\Gamma
    q \left( r, s, u (r,s) \right) \cdot \boldsymbol \nu \,
    \beta (r) \d{r}
  \end{displaymath}
  uniformly in $s \in ]-\infty,0[$ and
  \begin{displaymath}
    \lim_n\int_\Gamma
    \eta \left(u_b (r)\right) \, \beta (r) \d{r}
    =
    \int_\Gamma
    \eta \left(u_b (r)\right) \, \beta (r) \d{r}.
  \end{displaymath}
  Due to~\eqref{eq:7.39}, there exists a set $E_w$ of measure zero
  such that for all $\beta \in J$ there exists
  $\lim_{\substack{s \to 0^- \\s \notin E_w}} \int_\Gamma q \left( r,
    s, u (r,s) \right) \cdot \boldsymbol \nu \, \beta (r) \d{r}$
  and moreover
  \begin{displaymath}
    \lim_{\substack{s \to 0^- \\s \notin E_w}}
    \int_\Gamma
    q \left( r, s, u (r,s) \right) \cdot \boldsymbol \nu \,
    \beta (r) \d{r}
    \geq
    - \norma{\partial_u f}_{\L\infty ([0,T] \times \Omega \times \mathcal{U}; \reali^N)}
    \int_\Gamma
    \eta \left(u_b (r)\right) \, \beta (r) \d{r}.
  \end{displaymath}
  Note that the set $E_w$ depends on $w$ because $\eta$ and $q$
  depend on $w$. The above result can be extended to functions
  $\beta \in \L1 (\Gamma;\reali_+)$, so that for all
  $w \in \mathbb{Q}$ the quantity
  \begin{displaymath}
    \lim_{\substack{s \to 0^- \\s \notin E_w}}
    \int_\Gamma
    Q \left( r, s, u (r,s), w \right) \cdot \boldsymbol \nu \,
    \beta (r) \d{r}
  \end{displaymath}
  exists and
  \begin{equation}
    \label{eq:7.40}\!\!\!
    \lim_{\substack{s \to 0^- \\s \notin E_w}}
    \int_\Gamma
    Q \left( r, s, u (r,s), w \right) \cdot \boldsymbol \nu \,
    \beta (r) \d{r}
    \geq
    - \norma{\partial_u f}_{\L\infty ([0,T] \times \Omega \times \mathcal{U}; \reali^N)}\!\!
    \int_\Gamma \!\!
    H \left(u_b (r), w \right) \beta (r) \d{r}.
  \end{equation}
  Let $v \in \L\infty (\Gamma;\reali)$ and
  $\beta \in \L1 (\Gamma;\reali_+)$ be given, and let $(v_n)$ be a
  sequence of simple functions with values in $\mathbb{Q}$ which
  converges to $v$ almost everywhere in
  $\Gamma$. Obviously,~\eqref{eq:7.40} holds for all
  $w=v_n$. Moreover,
  \begin{displaymath}
    \lim_{n}
    \int_\Gamma
    Q \left( r, s, u (r,s), v_n (r)\right) \cdot \boldsymbol \nu \,
    \beta (r) \d{r}
    =
    \int_\Gamma
    Q \left( r, s, u (r,s), v (r)\right) \cdot \boldsymbol \nu \,
    \beta (r) \d{r}
  \end{displaymath}
  uniformly in $s \in ]-\infty,0[$ and
  \begin{displaymath}
    \lim_{n}
    \int_\Gamma
    H \left(u_b (r), v_n (r) \right) \beta (r) \d{r}
    =
    \int_\Gamma
    H \left(u_b (r), v (r) \right) \beta (r) \d{r}.
  \end{displaymath}
  Hence, 
  the following inequality holds
  \begin{displaymath}
    \lim_{\substack{s \to 0^- \\s \notin E_w}}
    \int_\Gamma
    Q \left( r, s, u (r,s), v (r) \right) \cdot \boldsymbol \nu \,
    \beta (r) \d{r}
    \geq
    - \norma{\partial_u f}_{\L\infty ([0,T] \times \Omega \times \mathcal{U}; \reali^N)}
    \int_\Gamma
    H \left(u_b (r), v (r) \right) \beta (r) \d{r}.
  \end{displaymath}
  Choosing $v=u_b$ and recalling the properties of the boundary entropy
  $H$ (see Definition~\ref{def:beef}) conclude the proof.
\end{proofof}

\begin{proofof}{Lemma~\ref{lem:mvbc}}
  The proof follows immediately from Lemma~\ref{lem:rebc}. Indeed, for
  $k \in \reali$ and $n \in \naturali \setminus \{0\}$, using the
  notation introduced in~\eqref{eq:I}, define the maps
  \begin{align*}
    \Delta^k (u,w) = \
    & \min_{z \in \mathcal{I}[w,k]} \modulo{u-z}
    \\
    H_n^k (u,w) = \
    & \left(\left(\Delta^k (u,w)\right)^2 + \frac1{n^2}\right)^{1/2} - \frac1n
    \\
    Q_n^k (t,x,u,w) = \
    & \int_w^u \partial_1 H_n^k (z,w) \, \partial_u f (t,x,z) \d{z}.
  \end{align*}
  It can be easily proved that, for all $k \in \reali$, the sequence
  of boundary entropy--entropy flux pairs
  $(H_n^k(u,w), Q_n^k(t,x,u,w))$ converges uniformly to
  $(\Delta^k(u,w), \mathcal{F} (t,x,u,w,k))$ as $n$ goes to
  $+\infty$. Applying~\eqref{eq:rebc}, with $Q$ replaced by $Q_n^k$,
  yields the thesis in the limit $n \to +\infty$, for all
  $k \in \reali$ and
  $\beta \in \L1 (]0,T[ \times \partial\Omega;\reali_+)$.
\end{proofof}

\begin{proofof}{Theorem~\ref{thm:ebln}}\vspace{-4mm}
  \paragraph{An E--solution is a BLN--solution.}
  It is sufficient to consider the following sequence of (classical)
  entropies: for $k \in \reali$
  \begin{displaymath}
    \eta_n (z) =  \sqrt{(z-k)^2+\dfrac{1}{n}},
  \end{displaymath}
  the corresponding entropy fluxes $q_n$ being defined as in
  point~2.~of Definition~\ref{def:eef}. A standard limiting procedure
  allows to obtain, in the limit $n \to + \infty$,
  \begin{align*}
    \eta_n (z) \to \ & \modulo{z-k}
    \\
    q_n (t,x,z) \to \ & \sgn (z-k) \, \left(f (t,x,z) - f (t,x,k)\right),
  \end{align*}
  so that, in the limit $n \to +\infty$,~\eqref{eq:e} becomes
  \begin{align*}
    & \int_0^T \!\!\int_{\Omega} \!\left[ \modulo{u (t,x) - k}
      \, \partial_t \phi (t,x)
      + \sgn\left(u (t,x) - k \right)
      \left( f\left(t,x,u (t,x) \right) - f\left(t,x, k \right)
      \right)
      \cdot \nabla \phi (t,x)\right]\! \d{x} \d{t}
    \\
    & + \int_0^T \int_{\Omega} \sgn\left(u (t,x) - k \right)
      \left[
      F\left(t,x,u (t,x) \right) - \div f \left(t,x,k \right)
      \right] \, \phi (t,x) \d{x} \d{t}
    \\
    & + \int_0^T \int_{\Omega} \sgn\left(u (t,x) - k \right)
      \div \left(
      f \left(t,x,u (t,x) \right) - f \left(t,x, k \right)
      \right)
      \, \phi (t,x) \d{x} \d{t}
    \\
    & +
      \int_{\Omega} \modulo{u_o (x)-k} \phi (0,x) \d{x}
    \\
    &-
      \int_0^T\!\!\! \int_{\partial \Omega}\!
      \sgn\left(u_b (t,\xi) - k \right)
      \left( f \left(t,\xi,u_b (t,\xi) \right) - f(t,\xi,k)\right) \cdot
      \nu (\xi) \, \phi (t,\xi) \d\xi \d{t}
    \\
    & +
      \int_0^T \int_{\partial \Omega}
      \sgn\left(u_b (t,\xi) - k \right)
      \left(
      f\left(t,\xi, u_b(t,\xi)\right) - f (t,\xi,\tr u (t,\xi))\right)
      \cdot
      \nu (\xi) \, \phi (t,\xi) \d\xi \d{t} \geq 0.
  \end{align*}
  Combining in the inequality above the second line with the third and
  the fifth one with the sixth yields exactly~\eqref{eq:bln}.

  \paragraph{A BLN--solution is an E--solution.}
  Assume that
  $\norma{u}_{\L\infty ([0,T] \times \Omega;\reali)} \leq U$.  It is
  immediate to see that $u$ satisfies~\eqref{eq:e} with
  $\eta(u) = \alpha\modulo{u - k} + \beta$, for any $\alpha > 0$ and
  $k, \beta \in \reali$.  Moreover, if $u$ satisfies~\eqref{eq:e} for
  two distinct locally Lipschitz continuous pairs $(\eta_1,q_1)$ and
  $(\eta_2, q_2)$, then the same inequality~\eqref{eq:e} holds for $u$
  with $(\eta_1+\eta_2, q_1+q_2)$. It can be proved by induction that
  $u$ satisfies~\eqref{eq:e} for any pair $(\eta,q)$ with $\eta$
  piecewise linear and continuous on $[-U,U]$.  Furthermore, if $u$
  satisfies~\eqref{eq:e} for the continuous pairs $(\eta_n,q_n)$ and
  the $\eta_n$ converge uniformly to $\eta$ on $[-U,U]$, then $u$
  fulfils~\eqref{eq:e} also for the pair $(\eta,q)$, where $q$ is
  defined as in point~2.~of Definition~\ref{def:eef}. To conclude,
  since any convex entropy $\eta$ is the uniform limit on $[-U, U ]$
  of piecewise linear and continuous functions, we obtain the proof.
\end{proofof}

\begin{proofof}{Lemma~\ref{lem:ebc}}
  The proof follows the lines of that
  of~\cite[Proposition~2.3]{bordo}. Indeed, let
  $\Phi \in \Cc1 (]0,T[ \times \reali^N; \reali_+)$ and
  $\psi_h \in \Cc1 (\bar \Omega; [0,1])$, with $\psi_h (\xi) = 1$ for
  all $\xi \in \partial \Omega$, $\psi_h (x)=0$ for all $x \in \Omega$
  with $B (x,h) \subseteq \Omega$, and
  $\norma{\nabla \psi_h}_{\L\infty (\Omega;\reali^N)} \leq 2/h$.
  Write~\eqref{eq:e} with $\phi (t,x) = \Phi (t,x) \, \psi_h (x)$ and
  take the limit as $h \to 0$. For any entropy--entropy flux pair
  $(\eta,q)$, thanks to the Dominated Convergence
  Theorem and to~\cite[Lemma~A.4 and Lemma~A.6]{bordo}, we get
  \begin{align*}
    \int_0^T\int_\Omega
    q\left(t,\xi,\tr u (t,\xi)\right)
    \cdot \nu (\xi) \, \Phi (t,\xi) \d\xi \d{t}
    -
    \int_0^T\int_\Omega
    q\left(t,\xi, u_b (t,\xi)\right)
    \cdot \nu (\xi) \, \Phi (t,\xi) \d\xi \d{t} &
    \\
    +
    \int_0^T\int_\Omega
    \eta' \left(u_b (t,\xi)\right)
    \left(
    f\left(t,\xi,u_b (t,\xi)\right)
    -
    f\left(t,\xi,\tr u (t,\xi)\right)
    \right)
    \cdot \nu (\xi) \, \Phi (t,\xi) \d\xi \d{t} &
    \geq 0.
  \end{align*}
  Hence, for any entropy--entropy flux pair
  $(\eta,q)$, \eqref{eq:ebc} holds almost everywhere on
  $]0,T[ \times \partial \Omega$.
\end{proofof}

\begin{proofof}{Lemma~\ref{lem:blnbc}}
  Proving that~\eqref{eq:6} holds is done in the same way as in the
  proof of~\cite[Proposition~2.3]{bordo}.

  It is immediate to prove that~\eqref{eq:6} reduces
  to~\eqref{eq:blnbc} when
  $k \in \mathcal{I}[\tr u (t,\xi), u_b (t,\xi)]$.  On the other hand,
  assume that~\eqref{eq:blnbc} holds and consider the various
  possibilities.
  \begin{itemize}
  \item If $k \leq \min\left\{ \tr u (t,\xi), u_b (t,\xi) \right\}$ or
    $k \geq \max\left\{ \tr u (t,\xi), u_b (t,\xi) \right\}$: the
    quantity
    \begin{displaymath}
      \sgn\left(\tr u (t,\xi) - k\right) - \sgn\left(u_b (t,\xi)- k\right)
    \end{displaymath}
    is equal to $0$, so~\eqref{eq:6} clearly holds.

  \item if $\tr u (t,\xi) \leq k \leq u_b (t,\xi)$: \eqref{eq:blnbc}
    reads
    $- ( f\left(t,\xi,\tr u (t,\xi)\right) - f (t,\xi,k)) \cdot \nu
    (\xi) \geq 0$,
    while $\sgn\left(\tr u (t,\xi) - k\right)=-1$ and
    $\sgn\left(u_b (t,\xi)- k\right)=+1$, so that~\eqref{eq:6} clearly
    holds.

  \item if $u_b (t,\xi) \leq k \leq \tr u (t,\xi)$: \eqref{eq:blnbc}
    reads
    $( f\left(t,\xi,\tr u (t,\xi)\right) - f (t,\xi,k)) \cdot \nu
    (\xi) \geq 0$,
    while $\sgn\left(\tr u (t,\xi) - k\right)=+1$ and
    $\sgn\left(u_b (t,\xi)- k\right)=-1$, so that~\eqref{eq:6} clearly
    holds.
  \end{itemize}
  The proof is completed.
\end{proofof}

\begin{proofof}{Proposition~\ref{prop:equivbc}}\vspace{-4mm}
  \paragraph{~\ref{item:1} $\boldsymbol\Rightarrow$~\ref{item:2}.} It
  is proved in Lemma~\ref{lem:mvbc}.

  \paragraph{~\ref{item:2} $\boldsymbol\Rightarrow$~\ref{item:4}.}
  From~\eqref{eq:mvbc} it follows that, for any
  $\beta \in \L1 (]0,T[ \times \partial \Omega;\reali_+)$ and for any
  $k\in \reali$,
  \begin{align*}
    &\int_0^T \int_{\partial\Omega}
      \mathcal{F}\left(t, \xi, \tr u (t, \xi), u_b (t,\xi), k\right)
      \cdot \nu (\xi) \, \beta (t,\xi) \d\xi\d{t}
    \\
    & = \
      \esslim_{\rho \to 0^+}
      \int_0^T \int_{\partial\Omega}
      \mathcal{F}
      \left(t, \xi, u (t, \xi -\rho \, \nu (\xi) ), u_b (t,\xi), k
      \right)
      \cdot \nu (\xi) \, \beta (t,\xi) \d\xi\d{t}.
  \end{align*}
  Therefore, there is a set $E \subseteq ]0,T[ \times \partial \Omega$
  of zero measure such that, for all $k\in \reali$ and for all
  $(t,\xi) \in (]0,T[ \times \partial\Omega) \setminus E $
  \begin{displaymath}
    \mathcal{F}\left(t, \xi, \tr u (t, \xi), u_b (t,\xi), k\right)
    \cdot \nu (\xi) \geq 0.
  \end{displaymath}

  \paragraph{~\ref{item:4} $\boldsymbol\Rightarrow$~\ref{item:5}.} It
  follows immediately from the definition~\eqref{eq:2} of
  $\mathcal{F}$.

  \paragraph{~\ref{item:5} $\boldsymbol\Rightarrow$~\ref{item:6}.} For
  any entropy--entropy flux pair $(\eta,q)$ and for any
  $(t,\xi) \in ]0,T[ \times \partial \Omega$, it holds
  \begin{align*}
    q (t,\xi,z) = \
    &
      q \left(t,\xi,w\right) + \int_w^z
      \eta' (\lambda) \, \partial_u f (t,\xi,\lambda)  \d\lambda
    \\
    = \
    & q \left(t,\xi,w\right) +
      \eta' \left(w\right) \left(
      f (t,\xi,z) -  f \left(t,\xi,w\right)
      \right)
      +
      \int_w^z \eta'' (\lambda) \,
      \left(f (t,\xi,z) -  f (t,\xi,\lambda) \right)  \d\lambda.
  \end{align*}
  The above formula and~\eqref{eq:blnbc} imply~\eqref{eq:ebc}.

  \paragraph{~\ref{item:6} $\boldsymbol\Rightarrow$~\ref{item:3}.}
  It is sufficient to apply~\eqref{eq:ebc} to any entropy--entropy
  flux pair $(\eta,q)$ with
  \begin{displaymath}
    \eta' \left(u_b (t,\xi)\right) = 0
    \quad \mbox{ and } \quad
    q \left(t, \xi, u_b (t,\xi)\right)=0.
  \end{displaymath}

  \paragraph{~\ref{item:3} $\boldsymbol\Rightarrow$~\ref{item:1}.}
  For any boundary entropy--entropy flux pair $(H,Q)$ and
  $\beta \in \L1 (]0,T[ \times \partial\Omega;\reali_+)$ it holds
  \begin{align*}
    & \int_0^T \int_{\partial\Omega}
      Q \left(t,\xi, \tr u (t,\xi), u_b (t,\xi)\right)
      \cdot \nu (\xi) \, \beta (t,\xi) \d\xi\d{t}
    \\
    & = \
      \esslim_{\rho \to 0^+}
      \int_0^T \int_{\partial\Omega}
      Q \left(t,\xi, u (t,\xi - \rho \, \nu (\xi)), u_b (t,\xi)\right)
      \cdot \nu (\xi) \, \beta (t,\xi) \d\xi\d{t}.
  \end{align*}
  The right hand side above is clearly positive, due to
  Definition~\ref{def:beef} and~\eqref{eq:dubois},
  proving~\eqref{eq:rebc}.
\end{proofof}

\begin{proofof}{Theorem~\ref{thm:mvbln}}
  Thanks to Theorem~\ref{thm:remv} and Theorem~\ref{thm:ebln}, we know
  that the following relations hold
  \begin{align*}
    \mbox{RE--solutions} \, \Longleftrightarrow \ & \, \mbox{MV--solutions}
    &\,\,\, \mbox{and} & \quad &
    \mbox{E--solutions} \, \Longleftrightarrow \ & \, \mbox{BLN--solutions.}
  \end{align*}
  Therefore, we now prove that a MV--solution is a BLN--solution and that
  an E--solution is a RE--solution.

  \paragraph{A MV--solution is a BLN--solution.}
  \hspace{-2.5mm}Let $k\in \reali$ and
  $\phi \!\in \Cc1 (]-\infty,T[ \times \pint \Omega; \reali_+)$.
  Adding~\eqref{eq:mv} with '$+$' and~\eqref{eq:mv} with '$-$' yields
  the following inequality:
  \begin{align}
    \nonumber
    & \int_0^T \int_{\Omega} \modulo{u (t,x) - k}
      \, \partial_t \phi (t,x) \d{x} \d{t}
    \\
    \nonumber
    & + \int_0^T \int_{\Omega}  \sgn (u (t,x) -k) \;
      \left(f\left(t,x,u (t,x) \right) -f\left(t,x,k \right)\right)
      \cdot \nabla \phi (t,x) \d{x} \d{t}
    \\
    \label{eq:add}
    & + \int_0^T \int_{\Omega} \sgn (u (t,x) -k)
      \left[
      F\left(t,x,u (t,x) \right) - \div f \left(t,x,k \right)
      \right] \, \phi (t,x) \d{x} \d{t}
    \\
    \nonumber
    & +
      \int_{\Omega} \modulo{u_o (x) - k }\, \phi (0,x) \d{x}
      \geq 0.
  \end{align}
  Fix $h > 0$ and consider as a test function
  $\Phi_h (t,x) = \phi (t,x) \left(1 - \psi_h (x)\right)$, with
  $\phi \in \Cc1 (]-\infty,T[ \times \reali^N; \reali_+)$ and
  $\psi_h \in \Cc1 (\bar\Omega; [0,1])$, with $\psi_h (\xi) = 1$ for
  all $\xi \in \partial \Omega$, $\psi_h (x) = 0$ for all
  $x \in \Omega$ with $B (x, h) \subseteq \Omega$, and
  $\norma{\nabla \psi_h}_{\L\infty (\Omega;\reali^N)} \leq 2/h$.
  Note that
  $\lim_{h \to 0} \left(1-\psi_h (x)\right)= \caratt{\Omega} (x)$ and
  $\Phi_h \in \Cc1 (]-\infty,T[ \times \pint \Omega; \reali_+)$. Using
  $\Phi_h$ into~\eqref{eq:add} yields
  \begin{align*}
    & \int_0^T \int_{\Omega} \modulo{u (t,x) - k}
      \left(1-\psi_h (x)\right) \partial_t \phi (t,x)
      \d{x} \d{t}
    \\
    &
      + \int_0^T \int_{\Omega}
      \sgn (u (t,x) -k)
      \left(f\left(t,x,u (t,x) \right) -f\left(t,x,k \right)\right)
      \cdot \nabla \phi (t,x) \left(1 - \psi_h (x) \right) \d{x} \d{t}
    \\
    &
      - \int_0^T \int_{\Omega}
      \sgn (u (t,x) -k)
      \left(f\left(t,x,u (t,x) \right) -f\left(t,x,k \right)\right)
      \cdot \phi (t,x) \, \nabla \psi_h (x)\d{x} \d{t}
    \\
    & + \int_0^T \int_{\Omega} \sgn (u (t,x) -k)
      \left[
      F\left(t,x,u (t,x) \right) - \div f \left(t,x,k \right)
      \right] \, \phi (t,x) \left(1-\psi_h (x)\right) \d{x} \d{t}
    \\
    & + \int_{\Omega}
      \modulo{u_o (x)-k} \, \phi (0,x) \left(1-\psi_h (x)\right) \d{x}
      \geq 0.
  \end{align*}
  Let now $h$ tend to $0$. Thanks to~\cite[Lemma~A.4 and
  Lemma~A.6]{bordo} we obtain
  \begin{align}
    \nonumber
    & \int_0^T \int_{\Omega} \modulo{u (t,x) - k}
      \partial_t \phi (t,x)
      \d{x} \d{t}
    \\
    \nonumber
    &
      + \int_0^T \int_{\Omega}
      \sgn (u (t,x) -k)
      \left(f\left(t,x,u (t,x) \right) -f\left(t,x,k \right)\right)
      \cdot \nabla \phi (t,x)  \d{x} \d{t}
    \\
    \label{eq:addlimit}
    &
      - \int_0^T \int_{\partial \Omega}
      \sgn (\tr u (t,\xi) -k)
      \left(f\left(t,\xi,\tr u (t,\xi) \right) -f\left(t,\xi,k \right)\right)
      \cdot \nu(\xi) \, \phi (t,\xi)\d{\xi} \d{t}
    \\
    \nonumber
    & + \int_0^T \int_{\Omega} \sgn (u (t,x) -k)
      \left[
      F\left(t,x,u (t,x) \right) - \div f \left(t,x,k \right)
      \right] \, \phi (t,x)  \d{x} \d{t}
    \\
    \nonumber
    & + \int_{\Omega}
      \modulo{u_o (x)-k} \, \phi (0,x)  \d{x}
      \geq 0.
  \end{align}
  Consider in particular the third line above:
  \begin{equation}
    \label{eq:3}
    - \int_0^T \int_{\partial \Omega}
    \sgn (\tr u (t,\xi) -k)
    \left(f\left(t,\xi,\tr u (t,\xi) \right) -f\left(t,\xi,k \right)\right)
    \cdot \nu(\xi) \, \phi (t,\xi)\d{\xi} \d{t}.
  \end{equation}
  Since $u$ is a MV--solution to~\eqref{eq:1}, we can apply
  Lemma~\ref{lem:rebc}, or equivalently~\ref{lem:mvbc}. Moreover,
  $u \in (\L\infty \cap \mathcal{TR}^\infty)([0,T] \times \Omega;
  \reali)$
  and thus Proposition~\ref{prop:equivbc} and Lemma~\ref{lem:blnbc}
  hold. Therefore, thanks to~\eqref{eq:6} and the positivity of the
  test function $\phi$, we get
  \begin{displaymath}
    [\eqref{eq:3}] \leq
    - \int_0^T \int_{\partial \Omega}
    \sgn (u_b (t,\xi) -k)
    \left(f\left(t,\xi,\tr u (t,\xi) \right) -f\left(t,\xi,k \right)\right)
    \cdot \nu(\xi) \, \phi (t,\xi)\d{\xi} \d{t},
  \end{displaymath}
  which inserted into~\eqref{eq:addlimit} yields~\eqref{eq:bln},
  concluding the proof.

  \paragraph{An E--solution is a RE--solution.}
  Let $\phi \in \Cc1 (]-\infty,T[ \times \pint \Omega;\reali_+)$ and
  $k \in \reali$. For any boundary entropy--entropy flux pair $(H,Q)$,
  set for any $t \in[0,T] $, $x \in \Omega$ and $z \in \reali$
  \begin{align}
   \label{eq:9}
    \eta (z) = \ & H (z,k), & q (t,x,z) = \ & Q (t,x,z,k).
  \end{align}
  By Definition~\ref{def:beef} of boundary entropy--entropy flux pair,
  $(\eta,q)$  is  an  entropy--entropy   flux  pair  with  respect  to
  $f$. Notice moreover that $\eta  (k)=0$. Since $u$ is an E--solution
  to~\eqref{eq:1}, it satisfies~\eqref{eq:e} with  the above choice of
  the  test  function, which,  thanks  to~\eqref{eq:9},  now reads  as
  follows:
  \begin{equation}
    \label{eq:reinside}
    \begin{aligned}
    & \int_0^T \int_{\Omega} \left[ H\left(u (t,x), k\right)
      \, \partial_t \phi (t,x)
      + Q\left(t,x,u (t,x),k \right) \cdot \nabla \phi (t,x)\right] \d{x} \d{t}
    \\
    & + \int_0^T \int_{\Omega} \partial_1 H \left(u (t,x), k \right)
      \left[
      F\left(t,x,u (t,x) \right) - \div f \left(t,x,k \right)
      \right] \, \phi (t,x) \d{x} \d{t}
    \\
    & + \int_0^T \int_{\Omega} \div Q \left(t,x,u (t,x), k \right)
      \, \phi (t,x) \d{x} \d{t}
    \\
    & +
      \int_{\Omega} H\left(u_o (x), k\right) \phi (0,x) \d{x}
      \geq 0.
    \end{aligned}
  \end{equation}

  Apply now Lemma~\ref{lem:ebc} and Proposition~\ref{prop:equivbc}. In
  particular, \eqref{eq:rebc} holds for any boundary entropy--entropy
  flux pair $(H,Q)$ and for any
  $\beta \in \L1 (]0,T[ \times \partial \Omega;\reali_+)$. We now
  follow the lines of the second part of the proof
  of~\cite[Theorem~7.31]{MalekEtAlBook} in order to prove that $u$
  satisfies~\eqref{eq:re}. The idea is to show that every $u$ which
  is a solution inside the domain $\Omega$, that
  is~\eqref{eq:reinside} holds, and which satisfies the boundary
  condition in a suitable way, is indeed a RE--solution.

  Define the following maps: for $z,w \in \reali$
  \begin{displaymath}
    \tilde H (z,w) =
    \begin{cases}
      \eta (z) - \eta (w) & \mbox{ if } z \leq w \leq k,
      \\
      0 & \mbox{ if } w \leq z \leq k,
      \\
      \eta (z) & \mbox{ if } w \leq k \leq z,
      \\
      \eta (z) & \mbox{ if } z \leq k \leq w,
      \\
      0 & \mbox{ if } k \leq z \leq w,
      \\
      \eta (z) - \eta (w) & \mbox{ if } k \leq w \leq z,
    \end{cases}
  \end{displaymath}
  and, for $t \in [0,T]$, $x \in \overline \Omega$,
  \begin{displaymath}
    \tilde Q (t,x,z,w) =
    \begin{cases}
      q (t,x,z) - q (t,x,w) & \mbox{ if } z \leq w \leq k,
      \\
      0 & \mbox{ if } w \leq z \leq k,
      \\
      q (t,x,z) & \mbox{ if } w \leq k \leq z,
      \\
      q (t,x,z) & \mbox{ if } z \leq k \leq w,
      \\
      0 & \mbox{ if } k \leq z \leq w,
      \\
      q (t,x,z) - q (t,x,w) & \mbox{ if } k \leq w \leq z.
    \end{cases}
  \end{displaymath}
  It is easy to see that
  $(\tilde H, \tilde Q) \in \C0 (\reali^2;\reali) \times \C0 ([0,T]
  \times \overline \Omega\times \reali^2; \reali^N)$.
  Define, for $n \in \naturali \setminus\{0\}$,
  \begin{displaymath}
    H_n (z,w) =
    \left\{
      \begin{array}{ll}
        \left.
        \begin{array}{ll}
          \eta (z) - \eta \left(w-\frac1n\right) & \mbox{ if } z \leq w -\frac1n
          \\
          0 & \mbox{ if } w -\frac1n \leq z \leq k + \frac1n
          \\
          \eta (z) - \eta \left(k+\frac1n\right)  & \mbox{ if } k+\frac1n \leq z
        \end{array}
                                                    \right\} & \mbox{ and } w \leq k,
        \\[30pt]
        \left.
        \begin{array}{ll}
          \eta (z) - \eta\left(k-\frac1n\right)  & \mbox{ if } z \leq k-\frac1n
          \\
          0 & \mbox{ if } k-\frac1n \leq z \leq w+\frac1n
          \\
          \eta (z) - \eta \left(w+\frac1n\right) & \mbox{ if } w + \frac1n \leq z
        \end{array}
                                                   \right\} & \mbox{ and } k \leq w.
      \end{array}
    \right.
  \end{displaymath}
  Then, $(\tilde H, \tilde Q)$ can be locally uniformly approximated
  by $(\tilde H_n, \tilde Q_n)$, defined as follows
  \begin{align*}
    \tilde H_n (z,w) = \
    & \int_\reali H_n (\lambda, w) \, \rho_{1/n} (z-\lambda) \d\lambda,
    \\
    \tilde Q_n (t,x, z,w) = \
    &
      \int_w^z \partial_1 \tilde H_n (\lambda, w) \, \partial_u f (t,x,\lambda)
      \d\lambda,
  \end{align*}
  where $\rho_{1/n}$ is a smooth mollifier. The pair
  $(\tilde H_n, \tilde Q_n)$ is clearly a boundary entropy--entropy
  flux pair. Since~\eqref{eq:rebc} holds, we have, for all
  $\beta \in \L1 (]0,T[ \times \partial \Omega; \reali_+)$,
  \begin{displaymath}
    \esslim_{\rho \to 0^+} \int_0^T \int_{\partial\Omega}
    \tilde Q_n
    \left(t, \xi,  u \left(t,\xi - \rho\, \nu (\xi)\right), u_b (t,\xi)\right)
    \cdot \nu (\xi) \, \beta (t,\xi) \d\xi\d{t} \geq 0,
  \end{displaymath}
  which becomes, as $n \to + \infty$,
  \begin{equation}
    \label{eq:7}
    \int_0^T \int_{\partial\Omega}
    \tilde Q
    \left(t, \xi,  \tr u \left(t,\xi \right), u_b (t,\xi)\right)
    \cdot \nu (\xi) \, \beta (t,\xi) \d\xi\d{t} \geq 0,
  \end{equation}
  where we use the hypothesis that $u$ admits a trace at the
  boundary. Going through all the cases in the definition of
  $\tilde Q$ and exploiting the properties of $\eta$ yield
  \begin{displaymath}
    \modulo{\tilde Q (t,x,z,w) - q (t,x,z)}
    \leq
    \norma{\partial_u f}_{\L\infty ([0,T]\times\Omega\times \mathcal{U};\reali^N)} \, \eta (w),
  \end{displaymath}
  where $\mathcal{U}$ is the interval $\mathcal{U}=[-U,U]$ with
  $U=\norma{u}_{\L\infty ([0,T] \times \Omega;\reali)}$.  Therefore,
  by~\eqref{eq:7}, for all
  $\beta \in \C1 (]0,T[ \times \partial \Omega; \reali_+)$
  \begin{equation}
    \label{eq:8}
    \begin{aligned}
      &  \int_0^T \int_{\partial\Omega} q
      \left(t, \xi, \tr u \left(t,\xi \right)\right)
      \cdot \nu (\xi) \, \beta (t,\xi) \d\xi\d{t} \\ \geq \ & -
      \norma{\partial_u f}_{\L\infty
        ([0,T]\times\Omega\times\mathcal{U};\reali^N)} \int_0^T
      \int_{\partial \Omega} \eta \left(u_b (t,\xi)\right) \beta
      (t,\xi) \d\xi \d{t}.
    \end{aligned}
  \end{equation}
  Fix $h>0$. Consider~\eqref{eq:reinside} with the test function
  $\Phi_h (t,x) = \phi (t,x) \left(1 - \psi_h (x)\right)$, with
  $\psi_h$ as in the first part of the proof of this Theorem, so that
  we obtain
  \begin{align*}
    & \int_0^T \int_{\Omega} H \left(u (t,x), k
      \right)  \left(1 - \psi_h (x)\right) \, \partial_t \phi (t,x)\d{x}\d{t}
    \\
    & + \int_0^T \int_{\Omega}
      Q \left(t,x,u (t,x), k \right)
      \cdot \nabla \phi (t,x) \, \left(1 - \psi_h (x)\right) \d{x} \d{t}
    \\
    & - \int_0^T \int_{\Omega}
      Q \left(t,x,u (t,x), k \right)
      \cdot \phi (t,x) \, \nabla \psi_h (x) \d{x} \d{t}
    \\
    & + \int_0^T \int_{\Omega}
      \partial_1 H \left(u (t,x), k \right) \, \left[ F \left(t,x,u
      (t,x)\right) - \div f\left(t,x,u (t,x)\right) \right] \phi
      (t,x) \, \left(1 - \psi_h (x)\right)  \d{x} \d{t}
    \\
    & + \int_0^T \int_{\Omega} \div Q \left(t,x,u (t,x), k \right)
      \, \phi (t,x) \, \left(1 - \psi_h (x)\right) \d{x} \d{t}
    \\
    & + \int_{\Omega} H\left(u_o (x), k \right) \, \phi (0,x)  \,
      \left(1 - \psi_h (x)\right) \d{x}
      \geq 0.
  \end{align*}
  Let now $h$ tend to $0$: by~\cite[Lemma~A.4 and Lemma~A.6]{bordo} we get
  \begin{align*}
    & \int_0^T \int_{\Omega} H \left(u (t,x), k
      \right) \, \partial_t \phi (t,x)\d{x}\d{t}
    \\
    & + \int_0^T \int_{\Omega}
      Q \left(t,x,u (t,x), k \right)
      \cdot \nabla \phi (t,x) \d{x} \d{t}
    \\
    & - \int_0^T \int_{\Omega}
      Q \left(t,\xi, \tr u \left(t,\xi \right), k \right)
      \cdot \phi (t,\xi) \d{\xi} \d{t}
    \\
    & + \int_0^T \int_{\Omega}
      \partial_1 H \left(u (t,x), k \right) \, \left[ F \left(t,x,u
      (t,x)\right) - \div f\left(t,x,u (t,x)\right) \right] \phi
      (t,x) \d{x} \d{t}
    \\
    & + \int_0^T \int_{\Omega} \div Q \left(t,x,u (t,x), k \right)
      \, \phi (t,x)\d{x} \d{t}
    \\
    & + \int_{\Omega} H\left(u_o (x), k \right) \, \phi (0,x) \d{x}
      \geq 0.
  \end{align*}
  Thanks to the definition of $q (t,x,z) = Q (t,x,z,k)$ and
  to~\eqref{eq:8} we obtain~\eqref{eq:re}, concluding the proof.
\end{proofof}

\bigskip

\noindent\textbf{Acknowledgement:} The present work was supported by
the PRIN 2015 project \emph{Hyperbolic Systems of Conservation Laws
  and Fluid Dynamics: Analysis and Applications} and by the
INdAM--GNAMPA 2017 project \emph{Conservation Laws: from Theory to
  Technology}.

\small{

  \bibliography{def}

  \bibliographystyle{abbrv}

}

\end{document}